\documentclass[12pt]{amsart}
\usepackage{amsfonts,amsmath,amsthm,amscd,amssymb,latexsym,cite,verbatim,texdraw,floatflt,
caption2}
\RequirePackage{hyperref}

\usepackage[a4paper]{geometry}

\theoremstyle
{plain}

\begin{document}

\title{Extremal balleans}

\author{Igor Protasov}

\maketitle
\vskip 5pt

{\bf Abstract.}
A ballean (or coarse space) is a set endowed with a coarse structure.
A ballean $X$  is called  normal if any two asymptotically  disjoint subsets of $X$  are  asymptotically separated.
We  say that a ballean $X$ is ultranormal (extremely normal) if any two unbounded subsets of $X$  are not asymptotically disjoint (every unbounded subset of $X$ is large).
Every maximal ballean is extremely normal and every extremely normal ballean is ultranormal,  but the converse statements do not hold.
A normal ballean is ultranormal   if and only if the Higson$^{\prime}$s corona of $X$  is a singleton.
A discrete ballean $X$ is ultranormal   if and only if $X$  is maximal.
We construct a series of concrete balleans with extremal properties.

\vskip 10pt

\vskip 5pt


\normalsize

\section{Introduction}

Let $X$  be a set. A family $\mathcal{E}$ of subsets of $X\times X$ is called a {\it coarse structure } if
\vskip 7pt

\begin{itemize}
\item{}   each $E\in \mathcal{E}$  contains the diagonal  $\bigtriangleup _{X}$,
$\bigtriangleup _{X}= \{(x,x): x\in X\}$;
\vskip 5pt

\item{}  if  $E$, $E^{\prime} \in \mathcal{E}$ then $E\circ E^{\prime}\in\mathcal{E}$ and
$E^{-1}\in \mathcal{E}$,   where    $E\circ E^{\prime}=\{(x,y): \exists z((x,z) \in  E,  \   \ (z, y)\in E^{\prime})\}$,   $E^{-1}=\{(y,x): (x,y)\in E\}$;
\vskip 5pt

\item{} if $E\in\mathcal{E}$ and $\bigtriangleup_{X}\subseteq E^{\prime}\subseteq E  $   then
$E^{\prime}\in \mathcal{E}$;
\vskip 5pt

\item{}  for any   $x,y\in X$, there exists $E\in \mathcal{E}$   such that $(x,y)\in E$.

\end{itemize}
\vskip 7pt

A subset $\mathcal{E}^{\prime} \subseteq \mathcal{E}$  is called a
{\it base} for $\mathcal{E}$  if, for every $E\in \mathcal{E}$, there exists
  $E^{\prime}\in \mathcal{E}^{\prime}$  such  that
  $E\subseteq E ^{\prime}$.
For $x\in X$,  $A\subseteq  X$  and
$E\in \mathcal{E}$, we denote
$E[x]= \{y\in X: (x,y) \in E\}$,
 $E [A] = \cup_{a\in A}   \   \   E[a]$
 and say that  $E[x]$
  and $E[A]$
   are {\it balls of radius $E$
   around} $x$  and $A$.

The pair $(X,\mathcal{E})$ is called a {\it coarse space} \cite{b14}  or a ballean \cite{b10}, \cite{b12}.

Each subset $Y\subseteq X$  defines the {\it subballean}
$(Y, \mathcal{E}_{Y})$, where $\mathcal{E}_{Y}$  is the restriction of $\mathcal{E}$  to $Y\times Y$.
A subset $Y$  is called {\it bounded} if $Y\subseteq E[x]$  for some $x\in X$  and
$E\in \mathcal{E}$.
\vskip 7pt

Given a ballean $(X, \mathcal{E})$, a subset $Y$  of $X$  is called \vskip 7pt

\begin{itemize}
\item{} {\it large} if  there exists $E\in \mathcal{E}$  such  that  $X= E [A]$;\vskip 5pt

\item{} {\it small} if  $(X\setminus A) \cap L$ is  large for each large subset $L$;\vskip 5pt

\item{} {\it thick} if,  for every  $E\in\mathcal{E}$,
there exists
  $a\in A$  such that $E[a]\subseteq A$. \vskip 5pt
\end{itemize}

\vskip 7pt

Every metric $d$ on a set $X$ defines the {\it metric ballean }  $(X, \mathcal{E} _{d})$,
 where $\mathcal{E} _{d}$  has the base $\{\{(x,y): d(x,y)\leq r  \} :  r\in\mathbb{R}^{+}\}$.
A ballean $(X, \mathcal{E})$  is called {\it metrizable}  if there exists a
 metric $d$  on $X $  such that  $\mathcal{E}= \mathcal{E} _{d}$.
 A ballean $(X, \mathcal{E})$ is metrizable if and only if $\mathcal{E}$ has a countable base [12, Theorem 2.1.1].
  Let $( X, \mathcal{E})$ be a ballean.
A subset $U$  of $X$ is called an {\it asymptotic neighbourhood}  of a subset $Y \subseteq X$   if, for every $E\in \mathcal{E}$,  $E[Y]\setminus  U$  is bounded.
\vskip 7pt

Two subsets $Y,  Z$  of $X$  are called\vskip 5pt

\begin{itemize}
\item{} {\it asymptotically disjoint} if, for every $E\in \mathcal{E}$,  $E[Y]\cap   E[Z]$  is bounded;\vskip 5pt

\item{} {\it asymptotically separated} if $Y, Z$  have  disjoint asymptotic neighbourhoods.

\end{itemize}

\vskip 7pt

A ballean $(X, \mathcal{E})$ is called {\it normal} [7]  if  any two
asymptotically disjoint subsets are asymptotically separated.
Every ballean $(X, \mathcal{E})$ with linearly  ordered  base of $\mathcal{E}$,  in particular, every metrizable ballean  is normal [7, Proposition 1.1].

A function  $f: X\longrightarrow \mathbb{R}$  is called {\it slowly oscillating} if,
 for any $E \in \mathcal{E}$  and $\varepsilon > 0$,  there exists a bounded subset $B$
 of $X$  such  that {\it diam }  $f(E[x])< \varepsilon$   for each $x\in  X  \  B$.
By [7, Theorem 2.2],  a ballean $(X, \mathcal{E}) $  is normal if and only if, for any two   disjoint  and asymptotically disjoint subsets $Y, Z$  of  $X$,   there exists a slowly oscillating function  $f: X \longrightarrow [0,1]$ such that
  $f | _{Y} = 0$  and  $f | _{Z} = 1$.

We say that an unbounded ballean $(X, \mathcal{E})$ is

\begin{itemize}
\item{} {\it  ultranormal } [1] if any two unbounded subsets of $X$  are not asymptotically disjoint;
\vskip 7pt

\item{} {\it  extremely normal} if any unbounded subset of $X$ is large.
\end{itemize}

\vskip 10pt

An unbounded ballean $(X, \mathcal{E})$,  is called {\it  maximal}  if $X$  is bounded in any stronger coarse structure.
By [13],  every unbounded subset of a maximal ballean is large and every small subset is bounded. Hence, every maximal ballean is extremely normal.

A family $\mathcal{B}$ of subsets of a set $X$, closed under finite  unions  and subsets,
 is called a {\it bornology} if  $\cup\mathcal{B}=X$.
For every ballean $X$,  the family $\mathcal{B}_{X}$
  of all bounded subsets of $X$ is a bornology.

A  ballean $(X, \mathcal{E})$ is called {\it discrete} (or {\it  pseudodiscrete}   [12],  or {\it thin}  [5])   if, for every $E \in \mathcal{E}$,   there exists
  $B\in \mathcal{B}_{X}$ such that  $E [x]=\{ x\}$  for each $x\in X \backslash B$.
Every bornology $\mathcal{B}$  on a set  $X$  defines the  discrete   ballean
 $(X, \{E_{B}:  B \in\mathcal{B}\})$, where
 $E_{B}[x] = B $
 if  $x\in B$   and
  $E_{B}[x] = \{x\} $ if $x\in X\backslash B$.
It follows that every discrete ballean  $X$  is uniquely determined by its bornology
$\mathcal{B}_{X}$,
 see [12, Theorem 3.2.1].
Every discrete ballean is normal  [12, Example 4.2.2].

For a  ballean $X,$  the following conditions are equivalent: $X$ is discrete, every function  $f: X\longrightarrow \{0,1\}$ is slowly  oscillating, every unbounded subset of $X$ is thick, see [12, Theorem 3.3.1] and  [3, Theorem 2.2].

An unbounded discrete ballean is called {\it ultradiscrete }  if the  family
$\{ X \backslash B:  B \in \mathcal{B}_{X}\}
$ is an ultrafilter  on $X$.
Every ultradiscrete  ballean is maximal [12, Example 10.1.2].
Thus, we have got
\vskip 7pt

ultradiscrete $\Longrightarrow$ maximal $\Longrightarrow$  extremely normal $\Longrightarrow$ ultranormal
\vskip 7pt

\noindent and no one arrow can be reversed,  see Section 3.

\section{Characterizations}

Let $(X, \mathcal{E})$ be an unbounded ballean.
We endow $X$  with the discrete topology, identify the
Stone-$\check{C}$ech compactification $\beta X$  of $X$ with the set of all ultrafilters on $X$ and denote
$X ^{\sharp}=\{p\in  \beta X: P$ is unbounded for each $P\in p\}$.
Given any $p, q \in X^{\sharp}$,  we write $p \parallel q$  if there exists
 $E\in \mathcal{E}$ such that
 $E[Q]\in p$ for each $Q\in q$.
By [7, Lemma 4.1],
$\parallel$
 is an equivalence relation on $X^{\sharp}$.
We denote by $\sim$
  the minimal (by inclusion) closed (in $X^{\sharp}\times  X^{\sharp} $) equivalence
  on $X^{\sharp}$  such that $\parallel \subseteq \sim$.
The compact Hausdorff space $X ^{\sharp} / \sim$  is called the {\it corona}  of  $(X, \mathcal{E})$, it is denoted by $\nu(X, \mathcal{E})$.
If every ball $\{ y\in X:  d(x,y) \leq  r\}$ in the metric space $(X, d)$
 is compact then $\nu(X, \mathcal{E }_{d})$ is
 the Higson's corona of $(X, d)$,  see [8] and [14].

We say that a function $f: X\longrightarrow  \mathbb{R}$  is {\it constant at infinity}  if there exists $c\in \mathbb{R}$  such that, for each  $\varepsilon > 0$, the set $\{ x\in X: | f(x) - c | > \varepsilon\}$  is bounded.
We say that $f$  is {\it almost constant} if $f |  _{X \backslash B}=$const  for some $B\in \mathcal{B} _{X}$.
If the bornology $\mathcal{B} _{X}$ is closed under countable unions then every function, constant at infinity, is almost constant.

If $f$  is constant at infinity then $f$ is slowly oscillating.
We denote so$(X, \mathcal{E})$  the set of all  bounded slowly oscillating functions on $X$.
For a bounded function  $f: X \longrightarrow \mathbb{R}$,  $f^{\beta}$   denotes the extension of $f$ to $\beta X$.

\vskip 10pt

{\bf Theorem 1}.  {\it For an unbounded normal  ballean  $(X, \mathcal{E})$, the following conditions are equivalent:
\vskip 7pt

$(1) \ $   every function $f \in so (X, \mathcal{E})$  is constant at infinity;

\vskip 7pt

$(2) \ $    $\nu(X, \mathcal{E})$  is a singleton;

\vskip 7pt

$(3) \  $    $(X, \mathcal{E})$ is ultranormal.
\vskip 10pt

Proof.}   $(3) \Longrightarrow  (1)  $.
We assume that a function $f \in \  so \ (X, \mathcal{E})$    is
 not constant at infinity.  Then there exists distinct $a, b \in \mathbb{R}$  and
 $p, q \in X^{\natural} $ such that
 $f^{\beta} (p)= a,  \  $
 $f^{\beta} (q)= b.  \  $
We put
 $\varepsilon= \frac{|a - b |}{4} $  and choose  $P\in p$,  $Q\in q$   such that
 $f(P)\subset (a-\varepsilon, a+ \varepsilon)$,  $f(Q)\subset (b-\varepsilon, b+ \varepsilon)$.
Given an arbitrary  $E \in \mathcal{E}$, we take $B\in \mathcal{B}_{X}$   such that
 $diam \  f ( E[x])< \varepsilon$  for each  $x\in X \setminus B$.
It follows that
$E[ P\setminus B]\cap E[ Q\setminus B] = \emptyset$
 so $P$  and $Q$  are asymptotically disjoint and we get a contradiction to (3).
\vskip 10pt

$(1) \Longrightarrow  (2)  $.
Let $p, q \in X^{\sharp}$.
By Proposition 8.1.4  from [12],
$p\sim q$  if and only if
$f^{\beta} ( p) = f^{\beta} ( q)$
for every $f\in  so \ (X, \mathcal{E})$.
\vskip 10pt

$(2)\Longrightarrow   (3) $.
We assume that  $(X, \mathcal{E}) $  is not ultranormal
and choose two disjoint and asymptotically disjoint
unbounded subsets $A, B$  of $X$.  We chose $p, q\in X^{\sharp} $  such
 that $A\in p$,  $B\in q$.
Since $X$ is normal, there exists $f\in  so \ (X, \mathcal{E})$  such that
$f\mid _{A}= 0 $,  $f\mid _{B}= 1 $.
Then
$f^{\beta} ( p) \neq f^{\beta} ( q)$
  and
  $\mid  \nu (X, \mathcal{E})   \mid > 1$
   by Proposition 8.1.4 from [12].
$ \ \ \Box$

\vskip 10pt

An unbounded ballean $X$ is called {\it  irresolvable}  [11]  if $X$  can not be
partitioned into two large subsets.
\vskip 10pt

{\bf Theorem 2}.  {\it
For an unbounded discrete ballean $(X, \mathcal{E})$, the following conditions are equivalent:
\vskip 10pt

$(1) \ $ $(X, \mathcal{E})$     is ultradiscrete; \vskip 7pt

$(2) \ $ $(X, \mathcal{E})$   is extremely normal; \vskip 7pt

$(3) \ $ $(X, \mathcal{E})$   is ultranormal; \vskip 7pt

$(4) \ $ $(X, \mathcal{E})$   is  maximal and  irresolvable.

\vskip 10pt

Proof.}
$(1) \Longrightarrow  (2)  $.
We denote by $p$  the ultrafilter $\{X\backslash B:  B \in \mathcal{B}_{X}\}$.
If $A$  is an unbounded subset of $X$  then $A\in p$  so $X\backslash  A$ is bounded and $A$  is large.

\vskip 5pt

$(2)\Longrightarrow   (3) $.
Let $A, C$  be unbounded subsets of $X$.
Since   $A$  is large, there exists $E\in \mathcal{E}$  such that $E[A] = X$  so $C\subseteq   E[A] $    and $A, C$  are not asymptotically disjoint.

\vskip 5pt

$(3)\Longrightarrow   (1) $.
 If $(X, \mathcal{E})$  is not ultradiscrete then  there exist two disjoint unbounded subsets $A, C$     of $X$.
Since $(X, \mathcal{E})$  is discrete,  $A$  and $C$  are asymptotically disjoint.

\vskip 5pt

$(1)\Longrightarrow   (4) $. This is  Theorem 10.4.5  from  [12].  $ \  \  \Box$

\vskip 15pt

 Let $\mathcal{B}$  be a bornology on a set  $X$.
Following  [1],  we say that a coarse  structure $\mathcal{E}$ is {\it compatible}   with  $\mathcal{B}$  if $\mathcal{B}$  is the bornology  of all  bounded  subsets of $(X, \mathcal{E})$.

Each bornology $\mathcal{B}$ on $X$  defines two coarse structures $\Downarrow \mathcal{B}$  and  $\Uparrow \mathcal{B}$,  the smallest and the largest  coarse  structures on $X$  compatible  with  $\mathcal{B}$.
Clearly,
$\Downarrow \mathcal{B}$
 is the discrete coarse structure  defined by $\mathcal{ B}$,  in  particular, $(X, \Downarrow \mathcal{B})$   is  normal.

The coarse structure
$\Uparrow \mathcal{B}$
 consists of all entourages
 $E\subseteq X\times X$  such that  $E= E ^{-1}$  and
   $E[B]\in \mathcal{B}$ for each  $B\in \mathcal{B}$.
But in  contrast  to $\Downarrow \mathcal{B}$,
   the coarse  structure
   $\Uparrow \mathcal{B}$
    needs not to be normal [2, Theorem 12].

If a ballean  $(X, \mathcal{E})$ is maximal then
.$\Uparrow \mathcal{B}_{X}=\mathcal{E}$.
It follows that every  maximal ballean  $(X, \mathcal{E})$ is uniquely determined by the bornology of bounded subsets:
  if $(X, \mathcal{E}^{\prime})$  is maximal and
  $\mathcal{B}_{(X, \mathcal{E})} = \mathcal{B}_{(X, \mathcal{E}^{\prime})} $
   then  $\mathcal{E}= \mathcal{E}^{\prime} $.

\vskip 10pt

{\bf Proposition 1}.  {\it
Let $(X, \mathcal{E})$  be an  extremely normal ballean and let $\mathcal{E}^{\prime}$  be coarse structure  on $X$  such that
$\mathcal{E}\subseteq \mathcal{E}^{\prime}$
 and
 $(X, \mathcal{E}^{\prime})$
 is maximal. Then
 $\mathcal{B}_{(X, \mathcal{E})} = \mathcal{B}_{(X, \mathcal{E}^{\prime})} $  and
 $ \mathcal{E}^{\prime} = \Uparrow \mathcal{B}_{(X, \mathcal{E})}$.

\vskip 7pt

Proof.}
We  assume the contrary and pick
$ A\in \mathcal{B}_{(X, \mathcal{E}^{\prime})}   \  \setminus \   \mathcal{B}_{(X, \mathcal{E}})$.
Since
$(X, \mathcal{E})$
is extremely normal,  $A$  is large in
$(X, \mathcal{E})$.
Hence, $A$  is large in $(X, \mathcal{E}^{\prime})$
 so $(X, \mathcal{E}^{\prime})$
  is bounded  and we get a contradiction  to maximality   of
  $(X, \mathcal{E}^{\prime})$.  $ \  \  \Box$
\vskip 10pt

Following [2], we say that a ballean $(X, \mathcal{E})  $  is {\it relatively maximal}  if  $\mathcal{E} = \Uparrow \mathcal{B} _{(X, \mathcal{E})}$.

We recall that two ultrafilters  $p, q\in \beta  X$  are  {\it incomparable} if,  for every $f : X\longrightarrow  X$,  we have $f^{\beta}(p )\neq  q$,  $f^{\beta}(q)\neq  p$.
\vskip 10pt

{\bf Proposition 2}.  {\it
 Let $(X, \mathcal{E})  $ be an unbounded discrete  ballean such that any two  distinct utrafilters from $X^{\sharp}$ are incomparable.
Then $(X, \mathcal{E})  $ is relatively maximal.
\vskip 10pt

Proof.}
We suppose  the contrary and choose a coarse structure $\mathcal{E}^{\prime}$  on $X$  such that $\mathcal{E}\subset \mathcal{E}^{\prime}$  and
 $\mathcal{B} _{(X, \mathcal{E})}= \mathcal{B} _{(X, \mathcal{E}^{\prime})}$.
Then there exists
$E^{\prime} \in \mathcal{E}^{\prime}$  such that the set  $Y= \{ y\in X: |E^{\prime}(y)|> 1\}$   is unbounded
in  $(X, \mathcal{E}^{\prime})  $.
For each $y\in Y$, we  pick
$f(y)\in E^{\prime}(y)$,   $f(y)\neq   y$  and extend  $f$  to $X$  by  $f(x)=x$  for each $x\in X\setminus  Y$.
We  take an arbitrary ultrafilter  $p\in (X, \mathcal{E})^{\sharp}$ such that $Y\in p$.
By the assumption,
 $f(P)$  is bounded in
 $(X, \mathcal{E})  $
  for some $P\in p$,  $P\subseteq Y$.
Then  $(E^{\prime})^{-1}  [f(P)]$ must be bounded in
$(X, \mathcal{E})  $
 and we get a contradiction with the  choice of $p$.  $  \   \   \   \Box$

\section{Constructions}

Given a family $\mathfrak{F}$  of subsets of $X\times X$,  we denote by $<\mathfrak{F}>$ the  intersections  of all  coarse  structures, containing each $F\cup \bigtriangleup_{X} $,  $F \in \mathfrak{F}$  and say that
$<\mathfrak{F}>$
 is {\it generated}  by $\mathfrak{F}$.
It is easy to see that $<\mathfrak{F}>$ has a base of subsets of the form $E_{0}\circ E_{1}\circ \ldots \circ E_{n}$, where
$$  F_{i}\in  \{ (F \cup \bigtriangleup_{X})\cup (F \cup \bigtriangleup_{X})^{-1} :  F \in \mathfrak{F}\}\cup
\{ (x,y) \cup \bigtriangleup_{X}: x,y \in X \} ,  \  n< \omega .$$
If $\mathcal{E} _{1}$ and $\mathcal{E} _{2}$ are coarse structures, we write
$\mathcal{E}_{1}  \bigvee \mathcal{E}_{2}$ in
place of $< \mathcal{E}_{1}  \cup \mathcal{E}_{2}>$ .
For lattices of  coarse structures, see [11].
\vskip 10pt

{\bf Proposition 3}.  {\it
Let $(X, \mathcal{E})  $
be an unbounded discrete ballean and let $p, q$  be distinct ultrafilters   from
$X^{\sharp}$.
If $(X, \mathcal{E})  $ is relatively maximal then $f^{\beta} (p)\neq  q$ for each
  $f: X\longrightarrow   X$ such that, for every $B\in \mathcal{B} _{(X, \mathcal{E})}$,
  $f(B)\in \mathcal{B} _{(X, \mathcal{E})}$
   and $f^{-1}(B)\in \mathcal{B} _{(X, \mathcal{E})}$.

\vskip 10pt

Proof.}
We   assume the  contrary and choose $f$  such that  $f^{\beta} (p)= q$.
  We put  $F=\{(x, f(x)):  x \in X\}$ and denote by $\mathcal{E}^{\prime}$ the coarse   structure generated by  $\mathcal{E}$ and $F$.
Since  $E$ is discrete and ${E}^{\prime}$ is not discrete, we have $\mathcal{E}\subset  {\mathcal{E}}^{\prime}$.

We take  $x\in X$  and   $E_{1}, \ldots , E_{n} \in \mathcal{E}\cup\{F\cup\bigtriangleup_{X}, \  \  F^{-1}\cup\bigtriangleup_{X}\}$.
Applying an induction  by $n$  and the assumptions
$f(B)\in \mathcal{B} _{(X, \mathcal{E})}$, 
$ \ f^{-1}(B)\in \mathcal{B} _{(X, \mathcal{E})}$
 for each
 $B \in \mathcal{B} _{(X, \mathcal{E})}$,
 we  conclude that
 $E _{1} \circ \ldots \circ E_{n} [x] \in \mathcal{B} _{(X, \mathcal{E})}$.
 Hence,
 $ \mathcal{B} _{(X, \mathcal{E})}= \mathcal{B} _{(X, \mathcal{E})^{\prime}}$
  and we  get  a contradiction to relative maximality of  $(X, \mathcal{E})$. $  \   \   \   \Box$
\vskip 10pt

{\bf Proposition 4}.  {\it
Every unbounded subballean   of maximal (extremely  normal,  ultranormal)  ballean is maximal (extremely normal, ultranormal).
\vskip 10pt

Proof.}
We prove only  the first  statement, the second and third are evident.

Let  $(X, \mathcal{E})$  be a maximal ballean, $Y$  be an unbounded subset of $X$.
We assume that
$(Y, \mathcal{E}_{Y})$
 is not maximal and choose a coarse structure  $\mathcal{E}^{\prime}$ on $Y$  such that
  $\mathcal{E}\mid_{Y} \subset  \mathcal{E}^{\prime}$  and
 $(Y, \mathcal{E}^{\prime})$
  is not bounded.
We put $\mathfrak{F} = <\mathcal{E}\cup \mathcal{E}^{\prime}>$.
Since $(X, \mathcal{E})$  is maximal and
$\mathcal{E}\subset \mathfrak{F}$, $(X,  \mathfrak{F})$ must be bounded.
On the other hand,  each bounded subset $B$ of
$(Y, \mathcal{E}^{\prime})$
 is bounded in
 $(X, \mathcal{E})$
  because otherwise $B$  is large in
  $(X, \mathcal{E})$
    so $B$  is large in
    $(Y, \mathcal{E}_{Y})$.
Now let $x\in X$,  and
$E_{1}, \ldots , E_{n}\in  \mathcal{E} \cup \{ E^{\prime} \cup \bigtriangleup_{X}: E^{\prime}\in \mathcal{E}^{\prime} \}.$
On induction by $n$, we see that
$E_{1} \circ \ldots\circ  E_{n} [x]$
 is bounded in $(X, \mathcal{E})$.
Hence $(X, \mathfrak{F})$ is not bounded and  we get a contradiction.  $ \  \  \Box$

\vskip 15pt
{\bf Proposition 5}.  {\it
Let  $\mathcal{E}, \mathcal{E}^{\prime}$  be coarse structures on a set $X$  such that  $\mathcal{E}\subseteq \mathcal{E}^{\prime}$ .
Then the following statements hold: \vskip 10pt

$(1)  \  \  $   if $(X, \mathcal{E})$ is extremely normal then $(X, \mathcal{E}^{\prime})$
 is extremely normal;
 \vskip 7pt

$(2)  \  \  $   if $(X, \mathcal{E})$ is
ultranormal and  $\mathcal{B}_{(X, \mathcal{E})} = \mathcal{B}_{(X, \mathcal{E}^{\prime})}$ then
$(X, \mathcal{E}^{\prime})$
 is ultranormal;
 \vskip 7pt

Proof.} $(1)  \  \  $
Let  $A$  be a subset of $X$.
If $A$  is unbounded in
$(X, \mathcal{E}^{\prime})$
 then $A$  is unbounded  in
 $(X, \mathcal{E})$.
If $A$  is large in $(X, \mathcal{E})$
 then $A$  is large in  $(X, \mathcal{E}^{\prime})$.
\vskip 7pt

$(2)  \  \  $   We assume that some unbounded  subsets $A, B$  of
$(X, \mathcal{E}^{\prime})$
 are asymptotically disjoint in
 $(X, \mathcal{E})$.
Then
$E^{\prime} (A)\backslash B \in \mathcal{B} _{(X, \mathcal{E}^{\prime})}$
 for each  $E^{\prime}\in \mathcal{E}^{\prime}$.
Since $\mathcal{E}\subseteq \mathcal{E}^{\prime}$
and
$\mathcal{B} _{(X, \mathcal{E}} = \mathcal{B} _{(X, \mathcal{E}^{\prime})}$,
we have
$E(A)\backslash B \in \mathcal{B} _{(X, \mathcal{E}^{\prime})}$
 for each $E\in \mathcal{E}$
 so  $A, B$  are asymptotically disjoint in $(X, \mathcal{E})$. $ \  \  \Box$


\vskip 15pt
{\bf Example  1}.
For every infinite  regular  cardinal $\kappa$,    we construct a  coarse  structure
$ \mathfrak{M}_{\kappa} $
  on  $\kappa$ such that  $(\kappa,  \mathfrak{M}_{\kappa})$
   is maximal and  $\mathcal{B} _{(\kappa,  \mathfrak{M}_{\kappa} )} = [\kappa]^{<\kappa}$.
We denote by $\mathfrak{F}$  the family of all coverings of
$\kappa$ defined by the rule:
$\mathcal{P} \in \mathfrak{F}$  if and  only if,  for each  $P\in \mathcal{P}$  and
 $x\in \kappa$,   $|P|< \kappa$  and
 $|\cup\{ P^{\prime}: x\in P^{\prime} $,  $P^{\prime}\in \mathcal{P}\}|< \kappa$.
Then $ \mathfrak{M}_{\kappa} $ is defined by the base
$\{ M _{\mathcal{P}}:  \mathcal{P}\in\mathfrak{F} \}$,
 where
 $M_{\mathcal{P}}= \{ (x,y): x\in P,  y\in P  \  $
 for some  $P\in\mathcal{P}\}$.
For general construction of  coarse  structures  by means of coverings, see [9] or [12, Section 7.5].
Clearly,
$\mathcal{B}_{(\kappa, \mathfrak{M_{\kappa}})} = [\kappa]^{<\kappa}$
and
$(\kappa, \mathfrak{M }_{\kappa})$
 is maximal  [12, Example 10.2.1].

\vskip 15pt

Let  $G$  be a group and let $X$ be a $G$-space with the action  $G\times X\longrightarrow  X  $,  $  \  (g, x)\longmapsto gx$.
A  bornology  $\mathcal{I}$  on $G$ is called a {\it group  bornology}  if,  for any  $A, B \in \mathcal{I}$,  we have  $AB \in \mathcal{I}$,  $ \ A ^{-1} \in \mathcal{I}$.
Every group bornology $\mathcal{I}$   on $G$  defines a coarse structure  $\mathcal{E} _{\mathcal{I}}$  on $X$  with  the base
$\{ E _{A}:   A\in \mathcal{I}, \  \  e \in A\}$,  where   $e$  is the identity of $G$,
 $E _{A} = \{ (x, y): y \in Ax\}$.
Moreover, every coarse structure on $X$  can be  defined in this way [6].

\vskip 15pt
{\bf Example  2}.
We define a coarse structure  $\mathcal{E}$  on $\omega$  such that the ballean  $( \omega, \mathcal{E})$ is  extremely normal but $(\omega, \mathcal{E})$
  is not maximal.
   Let  $S_{\omega}$  denotes the group  of all  permutations of $\omega$,  $\mathcal{I}= [S_{\omega}]^{<\omega}$,  $ \ \ \mathcal{E}= \mathcal{E}_{\mathcal{I}}$.
We show that every infinite subset $A$ of $\omega$ is large.
We partition $A$  and  $\omega$ into two infinite subset $A= A_{1} \cup A_{2}$,
$ \  \  W= W_{1} \cup W_{2}$
and choose two permutations  $f_{1},  f_{2}  $  of $\omega$  so that
\vskip 7pt

 $f_{1}(A_{1})= W_{1}$,  $ \  \  f_{1}(W_{1})=A_{1}$,   $ \  \   f_{1}(x) = x  \  \  $   for each  $  \  \  x \in \omega \setminus (A_{1}\cup W_{1})$,
\vskip 7pt

 $f_{2}(A_{2})= W_{2}$,  $ \  \  f_{2}(W_{2})=A_{2}$,   $  \  \  f_{2}(x) = x  \  \  $   for each  $  \  \  x \in \omega \setminus (A_{2}\cup W_{2})$
\vskip 10pt

Then we put $F= \{f_{1}, f_{2}, id\}$,  where $id$ is the identity permutation.
Clearly, $E_{F}[A]= \omega$  so $\mathcal{E}$ is  extremely normal.
To see that $( \omega, \mathcal{E})$  is not maximal, we note that
$\mathcal{E} \subseteq \mathfrak{M}_{\omega}$
  and choose a partition
  $\mathcal{P}= \{P_{n}: n\in\omega \}$
   of $\omega$ such that $|P_{n}|=n$.
Then, for each $n\in \mathbb{N}$,  there exist  $x\in \omega$  such that $M _{\mathcal{P} }[x]= n$.
For each
$H\in [S_{\omega}]^{< \omega}$
 and
 $x\in \omega$,
 we have $|E _{H} [x]| \leq |H|$.
Hence, $\mathcal{E}\subset \mathfrak{M}_{<\omega}$.

\vskip 15pt
{\bf Example  3}. Let $\kappa$  be a cardinal, $\kappa> \omega$.
We  construct  two coarse structures $\mathcal{E},  \mathcal{E}^{\prime}$  on $\kappa$ such that
$\mathcal{E}\subset \mathcal{E}^{\prime}$, $(\kappa,  \mathcal{E})$   is ultranormal  but not extremely normal,
$(\kappa,  \mathcal{E}^{\prime})$
 is extremely  normal but not maximal. Let
 $S_{\kappa}$
  denotes the group of all  permutations
$\kappa$,   $\mathcal{I}= [S_{\kappa}]^{<\omega}$,
$\mathcal{E}= \mathcal{E}_{\mathcal{I}}$.
Clearly, $\mathcal{B}_{(\kappa, \mathcal{E})}= [\kappa]^{<\omega}$.
Let $A, B$ be infinite subsets of $\kappa$.
We take countable subset
 $A^{\prime}\subseteq A$,  $B^{\prime}\subseteq B$  and use  argument from Example 2 to choose
 $F=\{ f_{1}, f_{2}, id \}$ so that
 $B^{\prime}\subseteq E_{F} [A^{\prime}]$.
It follows that  $A, B$ are not asymptotically  disjoint in
$(\kappa, \mathcal{E})$
  so
  $(\kappa, \mathcal{E})$
   is ultranormal.
If $A$ is a countable subset of $\kappa$ then
$| E_{H} [A]|= \omega$
  for each
  $H\in [S_{\kappa}]^{<\omega}$
   so
   $(\kappa, \mathcal{E})$
   is not extremely normal.

We denote by $\mathcal{F}$ the discrete coarse  structure on $\kappa$ defined by the
bornology
$[\kappa]^{<\kappa}$
 and put
 $\mathcal{E}^{\prime}= \mathcal{F}\vee \mathcal{E}$.
If $A$  is an unbounded subset of  $\mathcal{E}^{\prime}$  then
$|A|= \kappa$.  Applying the arguments from Example 2, we see that
$(\kappa, \mathcal{E}^{\prime})$
  is extremely normal. The ballean
  $(\kappa, \mathcal{E}^{\prime})$
    is not maximal because
    $\mathcal{E}^{\prime}\subset\mathcal{F}\vee\mathcal{E}_{J}$,
       where $J= [S_{\kappa}]^{\leq\omega}$.

 \section{Comments}

1.	Let  $(X, \mathcal{E})$  be a ballean,    $A$ and $B$  be subsets  of $X$.
Following [9],  we write  $A \delta B$  if and only if there exists
 $E \in \mathcal{E} $ such that  $A\subseteq  E [B]$,  $B\subseteq  E[A]$.

Let  $\mathcal{E}, \mathcal{E}^{\prime}$ be coarse  structures on a set $X$  such that
$\mathcal{B}_{(X,\mathcal{E})}= \mathcal{B}_{(X,\mathcal{E}^{\prime})} $.
If  $\mathcal{E}$  and $\mathcal{E}^{\prime}$ have  linearly ordered bases and  either
$so (X, \mathcal{E})= so (X, \mathcal{E}^{\prime})$   or
$\delta_{(X,\mathcal{E})}= \delta_{(X,\mathcal{E}^{\prime})} $
then
$\mathcal{E}= \mathcal{E}^{\prime}$,
 see [4, Theorem 2.1]  and  [3, Theorem 4.2].

We take the coarse structures $\mathcal{E}$  and $\mathfrak{M}_{\omega}$ on $\omega$ from Example 2.
By Theorem 1,  so
$so (\omega, \mathcal{E})= so (\omega, \mathfrak{M}_{\omega})$.
Since
$\mathcal{B} _{(\omega, \mathcal{E})}= \mathcal{B} (_{\omega, \mathfrak{M}_{\omega})} = [\omega]^{<\omega}$
 and
 $ (\omega, \mathcal{E}), (\omega, \mathfrak{M}_{\omega})$
 are  extremely normal, we have
 $\delta _{(\omega, \mathcal{E})}= \delta _{(\omega, \mathfrak{M}_{\omega})}$.
By the construction,
$\mathcal{E}\neq \mathfrak{M}_{\omega}$.

In this connection we remind  Question  7.5.1 from [12]: does
$\delta _{(X, \mathcal{E})}= \delta _{(X, \mathcal{E}^{\prime})}$
 imply
 $so (X,  \mathcal{E})= so (X,  \mathcal{E}^{\prime})$ ,
 and give the affirmative answer to this question.
 \vskip 10pt

Clearly,
$\delta _{(X, \mathcal{E})}= \delta _{(X, \mathcal{E}^{\prime})}$
 implies
 $\mathcal{B} _{(X, \mathcal{E})}= \mathcal{B} _{(X, \mathcal{E}^{\prime})}$.
Assume that there exists
$f\in so (X, \mathcal{E})\setminus so  (X, \mathcal{E}^{\prime})$.
Then there exist $\varepsilon > 0$
 and  $E^{\prime}\in \mathcal{E}^{\prime}$
  such that, for each
  $B\in\mathcal{B}_{(X, \mathcal{E}^{\prime})}$,
   one  can find
   $y_{B}$,  $z_{B}\in X\setminus B$
    such that
    $(y_{B},  z_{B})\in E^{\prime}$
     but
     $|f(y_{B})-  f (z_{B})|>\varepsilon$
We put
$Y= \{ y_{B}: B\in \mathcal{B}_{(X, \mathcal{E}^{\prime})}  \}$
 and choose a function
 $h: X\longrightarrow X$
  such that
 $(y, h(y))\in E^{\prime}$,  $|f(y)- f(h (y)) |> \varepsilon$
  for each  $y\in Y$.
 Let $p\in X^{\sharp} $  and  $Y\in p$, $q$ be an  ultrafilter with  the base
 $\{h(P):  P\in p\}$.
Then
$|f^{\beta}(p)- f^{\beta} (q)|\geq \varepsilon$
 and there exists
  $P^{\prime}\in p$  such that
  $|f(x)  -  f(y)|> \frac{\varepsilon}{2}$
   for all
   $x\in P^{\prime}$,  $y\in h (P^{\prime})$.
    Since $f\in  so  (X, \mathcal{E})$,  $P^{\prime}$  and  $h(P^{\prime}) $  are not  close in
    $(X, \mathcal{E})$  but $P^{\prime}$,  $h(P^{\prime}) $
     are  close in $(X, \mathcal{E}^{\prime})$ .

\vskip 15pt

2. Every  group $G$  has  the  natural {\it finitary coarse  structure}  $\mathcal{E }_{fin}$  with the base
$\{E _{F}: F\in  [G] ^{< \omega} , \  e\in F\}$,  where
$E_{F} = \{(x, y): y\in F x\}$.
Let $G$  be an uncountable abelian group.
By [4,  Corollary  3.2],
$(G, \mathcal{E}_{fin})$
 is not normal but every  function
 $f\in  so (G, \mathcal{E} _{fin})$ is constant at infinity.  This example shows  that the assumption of  normality in  Theorem 1 can not be omitted.

\vskip 15pt

3. Example 3 shows that Proposition 1 does not hold for  ultranormal  ballean in place of extremely normal.
Indeed, $(\kappa , \mathcal{E})$ is ultranormal, $\mathcal{B}_{(\kappa , \mathcal{E})}= [\kappa]^{<\omega} $,
$\mathcal{E}\subset \mathcal{E}^{\prime}$, $\mathcal{B}_{(\kappa , \mathcal{E}^{\prime})}= [\kappa]^{<\kappa} $.
We take a maximal coarse structure  $\mathcal{E}^{\prime\prime}$ such that $\mathcal{E}^{\prime}\subset \mathcal{E}^{\prime\prime}$.
Then $\mathcal{E}\subset \mathcal{E}^{\prime\prime}$.
 and
 $\mathcal{B}_{(\kappa , \mathcal{E})}$
  $\neq \mathcal{B}_{(\kappa , \mathcal{E}^{\prime\prime})}$.
\vskip 15pt

4. Following [1], we say that a ballean  $(X, \mathcal{E})$  {\it has  bounded  growth }  if there is a mapping
$f: X\longrightarrow \mathcal{B}_{X}$  such that
\vskip 7pt

\begin{itemize}
\item{}    $\cup_{x\in B} f(x)\in\mathcal{B}_{X}$ for each $B\in\mathcal{B}_{X}$;
\vskip 7pt

\item{}  for each   $E\in\mathcal{E}$,  there exists $C\in\mathcal{B}_{X}$ such that  $E[x]\subseteq f(x)$  for each
$x\in X\setminus C$.
\end{itemize}

\vskip 10pt

Clearly, every discrete  ballean has  bounded  growth $(f(x)= \{x\})$.
Let us take the maximal ballean $(\omega, \mathfrak{M}_{\omega})$
  from  Example 1. Since  each ball in $(\omega, \mathfrak{M}_{\omega})$
   is finite, one  can use  the  diagonal   process to  show  that $(\omega, \mathfrak{M}_{\omega})$
    is not of bounded growth.

\vspace{6 mm}

CONTACT INFORMATION

I.~Protasov: \\
Faculty of Computer Science and Cybernetics  \\
        Kyiv University  \\
         Academic Glushkov pr. 4d  \\
         03680 Kyiv, Ukraine \\ i.v.protasov@gmail.com

\end{document}